\def\d{\delta}

\def\la{\lambda}

\def \dd#1{{\bf#1}}



\def\ouv#1{\smash{\mathop{#1}\limits^{\lower 1pt\hbox
{$\scriptscriptstyle\circ$}}}}

\def\hfl#1#2{\smash{\mathop{\hbox to 12mm{\rightarrowfill}}
\limits^{\scriptstyle#1}_{\scriptstyle#2}}}


\long\def\eno#1#2{\par\smallskip{\bf{#1}}{\it\ {#2}}\par\medskip}

\def\tit#1{\vskip 5mm plus 1mm minus 2mm {\tir #1}
		\vskip 3mm plus 1mm minus 2mm}

\def\stit#1{\vskip 3mm plus 1mm minus 2mm {\bf{#1}}
		\smallskip}

\long\def\dem#1#2{{\bf {#1}}{\ {#2}$\diamondsuit$}\medskip}

\font\tir=cmbx10 at 12pt

\def\ref#1#2#3#4{{\bf #1}{\ #2}{\it ,\ #3}{,\ #4}\medskip}


\def \picture #1 by #2 (#3){\midinsert \centerline 
{\vbox to #2{\hrule width #1 heigth 0pt 
depth 0pt \null \vfill \special {picture #3}}}\endinsert }

\def\scaledpicture #1 by #2 (#3 scaled #4) {{
\dimen0 =#1 \dimen1 =$2
\divide \dimen0 by 1000 \multiply \dimen0 by #4
\divide \dimen1 by 1000 \multiply \dimen1 by #4
\picture \dimen0 by \dimen1 (#3 scaled $4)}}

\def\figure #1 #2 #3 {\midinsert \vglue 3mm 
{\vbox to #3 {\hrule width 6cm height 0cm depth 0cm \vfill
{\special {picture #1 scaled #2}}}}\vglue 2mm \endinsert}

\magnification=1200

\bigskip
\bigskip

{\centerline {\tir  RIGIDITY OF AMN VECTOR SPACES  }}
  
\bigskip

{\centerline {{\bf E. Mu\~noz Garcia} \footnote {(*)} 
{U.C.L.A. Department of Mathematics, 405, Hilgard Ave.
Los Angeles CA-90095-1555, U.S.A.; e-mail: munoz@math.ucla.edu}}}

\bigskip
\bigskip
\bigskip

{\bf Abstract:} {\it A metric vector space is asymptotically metrically
normable (AMN) if there exists a norm 
asymptotically isometric to the distance.
We prove that AMN vector spaces
are rigid in the class of metric vector spaces under asymptotically
isometric perturbations. This result follows from a 
general metric normability 
criterium.
 If the distance is translation invariant and 
satisfies an approximate multiplicative condition then there
exists a lipschitz equivalent norm. Furthermore, we give necessary and 
sufficient conditions for the distance to be asymptotically 
isometric to the norm.}

\bigskip
\bigskip


\bigskip
\bigskip


\bigskip
\bigskip

\tit {Introduction.}

The geometrization of algebraic structures is a fruitful
modern idea born from the collusion of synthetic
geometry  and classical algebra.
A successful example is Gromov's notion of hyperbolic groups 
that has shed 
new light on classical group theory (see [Gr], [GH]).
In this article we investigate metric vector spaces
from a metrical point of view.

\bigskip

Let $(E,d)$ be a metric space. 
The distance $d$ is {\it asymptotically isometric} to 
a distance $\d$ if for any $C_1>1$ there exists $C_2 \geq 0$ such 
that for any $x,y \in E$ we have
$$
C_1^{-1} \d (x,y) -C_2\leq d(x,y) \leq C_1 \d (x,y) +C_2.
$$
Two metric spaces $(E_1 , d_1)$ and $(E_2, d_2)$ are asymptotically 
isometric if there exists a one-to-one correspondence $\varphi : 
E_1 \to E_2$ such that $d_1$ is asymptotically isometric to the 
distance $\d =\varphi^* \ d_2$,
$$
\d (x,y) = d_2 (\varphi (x),\varphi (y)) \ .
$$

A metric vector space $(E,d)$ is a topological vector space
whose topology is generated by the distance $d$.

\eno {Definition 1.}{ A metric space $(E,d)$ is
asymptotically metrically normable (AMN) 
if there exists
a norm $|| \ ||$ on $E$ such that $(E,|| \ ||)$
is asymptotically isometric 
to $(E,d)$.}

 Our main theorem is:

\eno {Theorem 1.}{ AMN vector spaces are 
asymptotically isometrically rigid.
More precisely, let $(E,d)$ be a metric vector space.
 If $E$ is 
asymptotically isometric to an AMN vector space, then
$E$ is AMN.
 }

Theorem 1 is a consequence of a metric normability criterium for 
metric vector spaces (theorems 2 and 3.) 

\bigskip

It is well known that a Hausdorff topological vector 
space is metrizable if and only if the origin has a
countable neighborhood base. In this case there 
exists a translation invariant distance generating the 
topology (see for example [Sch] p.28.) 
A  Hausdorff topological vector space is normable if and only 
if the origin has a bounded convex neighborhood ([Sch] 
p.41.) Recall in a topological vector space a set $A$ is bounded if 
for each neighborhood $U$ of $0$ there exists a scalar $\lambda$ 
such that $A\subset \lambda U$.

These conditions are sharp but not always useful. For instance,
given a distance $d$ defining the topological vector space 
structure there is no effective way of determining the existence 
of a convex neighborhood of $0$. Also an asymptotically isometric 
perturbation does not preserve convex sets (even those at "infinity".)
The purpose of the following theorems is to exhibit explicit metric
conditions on the distance that imply the 
existence of a lipschitz equivalent norm.
Similar ideas were used by the author in the study of 
H\"older absolute values over a field (see [Mu].)

From now on we consider vector spaces over a locally compact
valued field $K$ of characteristic zero ($\dd Q \subset K$) 
and such that $(\dd Q , |.|)$ is archimedian. From the 
classification of locally compact fields (see for example 
[We] chapter I.3) we have  that $K$ is an $\dd R$-field, i.e. $K=\dd R$,
$K=\dd C$ or $K=\dd H$ the field of quaternions.
Note that the group of units $\dd U =\{ u\in K ; |u|=1\}$ is 
a compact topological group.

All results and proofs as given are valid for modules over 
a valuated ring $(A,|.|)$ with unit such that 
$A$ is of characteristic $0$ ($\dd Q \subset A$), $(A, |.|)$ 
is locally compact and the restriction of the absolute
value to $\dd Q$ is archimedian.

\medskip

\eno {Definition 2.}{Let $(E,d)$ be a metric vector space.
Given a constant $C_0 \geq 0$ the distance is $C_0$-translation
invariant if translations are $C_0$-isometries, that is
for all $x,y, z \in E$,
$$
d (x,y) -C_0 \leq d (x+z,y+z) \leq d (x , y) +C_0.
$$
Note that this is equivalent to the right hand side inequality, for 
all $ x,y, z \in E$, 
$$
d (x+z,y+z) \leq d (x , y) +C_0.
$$
}

\eno {Definition 3.} {
A distance $d$ on $E$ is lipschitz
equivalent (or $(C_1, C_2)$-lipschitz equivalent) 
to another distance $\d$ on $E$ if there exists  
$C_1 \geq 1$ and $C_2 \geq 0$ such that for $x,y\in E$,
$$
C_1^{-1} \d (x,y) -C_2 \leq d (x,y) \leq C_1 \d (x,y) +C_2
$$
}

\eno {Definition 4.}{A metric vector space $(E,d)$ is 
metrically normable (MN) if $d$ is lipschitz equivalent to 
a norm on $E$.}

\eno {Definition 5.}{Let $E$ be a $K$ vector space.
A distance $d$ on $E$ is lipschitz multiplicative
(or $(C_1, C_2 , C_3)$-lipschitz multiplicative) if there
are three constants, $C_1 \geq 1$, $C_2 \geq 0$ and 
$C_3 \geq 0$, such that for any $\lambda \in K$, $x,y \in E$,
we have
$$
C_1^{-1} |\la | d (x,y) -C_2 |\la | -C_3 \leq d (\la x , \la y ) 
\leq C_1 |\la | 
d (x,y)+C_2 |\la | +C_3.
$$
}

Notice that a MN vector space is lipschitz multiplicative. More precisely,
of $d$ is $(C_1, C_2)$-lipschitz equivalent to a norm, then 
$d$ is $(C_1^2 , C_1 C_2, C_2)$-lipschitz multiplicative.

We denote by $\mu$ the (rigth invariant) Haar measure on the compact 
group $(\dd U , .)$ normalized to have total mass $1$.

\eno {Theorem 2.}{Let $(E,d)$ be a   
metric vector space over $K$
with the distance $d$ $C_0$-translation 
invariant and $(C_1 ,C_2 , C_3)$-lipschitz
multiplicative.

Let $E_0$ be the maximal subspace of $E$ where the distance 
$d$ is bounded. Then $E_0$ is a closed subspace of $E$ and 
the quotient is a metrizable vector space. The Hausdorff
distance on classes modulo $E_0$ induced by $d$ defines a 
distance $D$ in the quotient. 

The vector space $(E/E_0 ,D)$
is metrically normable by a norm $||.||$
 which is $(C_1^2,C_2')$-lipschitz equivalent 
to $d$ with $C_2'=C_1C_2+C_1C_3+C_2$, i.e. for $x,y \in E$,
$$
C_1^{-2} \ d (x,y) -C_2'\leq ||x-y|| \leq C_1^2 \  d (x,y)+C_2'.
$$

Moreover, the norm $||.||$ can be defined by
$$
||\bar x ||=\lim_{n\to +\infty  } {1\over n}  \int_{\dd U} 
\ d (nux, 0) d\mu (u) .
$$

In particular, if the distance $d$ is unbounded in all 
non-trivial subspaces, then $E_0=\{ 0\}$ and $E$ is metrically
normable with a norm lipschitz equivalent to $d$.
}

As  mentioned before, we can construct in any metrizable	
topological vector space a translation invariant distance 
generating the topology thus the problem of metric normability 
is reduced by theorem 2 to construct such a distance that 
satisfies an approximate scalar multiplicative property and 
that is unbounded in non-trivial subspaces.

\bigskip
\bigskip

Using similar ideas we can characterize completely those 
distances that are asymptotically isometric to a norm.

\eno {Definition 6.}{Let $(E,d)$ be a  
vector space over $K$. The distance $d$ is asymptotically 
multiplicative if for any $C_1 >1$, there exists
$C_2 \geq 0$ and $C_3\geq 0$ such that for $x,y \in E$, 
we have
$$
C_1^{-1} \ |\la | \ \ d (x,y) -C_2 |\la |-C_3 \leq 
d (\la x , \la y) \leq C_1 \ 
|\la | \ d (x,y) +C_2 |\la | +C_3.
$$
}

Our last theorem gives  a necessary and sufficient 
condition for a distance to be asymptotically isometric 
to a norm.

\eno {Theorem 3.}{Let $(E,d)$ be a   
metric vector space over $K$.
The distance $d$ is asymptotically isometric to a norm $||.||$
if and only if $d$ is  asymptotically multiplicative, 
$d$ is unbounded in non-trivial subspaces and   
for any $C_1 >1$ there exists $C_0 \geq 0$ such that for any 
$n\geq 2$ and $x_1 , \ldots , x_n , y_1 , \ldots , y_n \in E$,  
we have
$$
d \left ( \sum_{i=1}^n x_i , \sum_{i=1}^n y_i \right )
\leq C_1 \sum_{i=1}^n d (x_i , y_i) + n C_0.
$$
In that case the norm $||.||$ can be obtained as 
described in theorem 2.}

It is not difficult to see that the conditions stated 
are necessary.
The last condition is related to the condition of 
translation invariance in the first theorem in the 
following way (more precisely see lemma 1 below):
All translations are isometries if and only  
if for any $x_1, x_2, y_1, y_2 \in E$,
$$
d (x_1+x_2, y_1+y_2) \leq d (x_1, y_1) +d (x_2, y_2).
$$

We first prove theorem 2, then theorem 3 follows along the 
same lines and finally theorem 1 follows from theorem 3.

\null
\vfill
\eject

\tit {1) Proof of theorem 2.}

For the first part we only need to assume that $d$ 
is $C_0$-translation invariant.

\eno {Proposition 1.}{Let $(E,d)$ be a metric vector 
space such that $d$ is $C_0$-translation invariant. 
Then for any $x,y \in E$ the 
following limit exists
$$
\d (x, y) =\lim_{n\to +\infty} {1\over n } 
d (nx , ny),
$$
and we have
$$
\d (x,y) \leq d (x,y)+2 C_0.
$$
}

\eno {Lemma 1.}{Let $(E,d)$ be a metric vector 
space. If the distance $d$ is 
$C_0$-translation invariant then we have
for any $x_1, x_2, y_1, y_2 \in E$,
$$
d (x_1+x_2, y_1+y_2) \leq d (x_1, y_1) +d (x_2, y_2)+2 C_0.
$$
Conversely, if we have the previous inequality then the 
distance $d$ is $2C_0$-translation invariant.
}

\dem {Proof Lemma 1.}{
We have 
$$
\eqalign {
d (x_1+x_2, y_1+y_2) &\leq d (x_1 , y_1 +y_2 -x_2)  +C_0\cr
&\leq d (x_1 ,y_1) +d (y_1 , y_1+y_2-x_2) +C_0\cr
&\leq d (x_1 ,y_1) + d (x_2 +(y_1-x_2),y_2 +(y_1-x_2)) +C_0\cr
&\leq d (x_1, y_1) +d (x_2, y_2) +2C_0\cr
}
$$ 
Conversely, the inequality with $x_2=y_2$ proves
that $d$ is $2C_0$-translation invariant.}

\dem {Proof of proposition 1.}{ Consider for $n\geq 0$,
$$
a_n =d (nx , ny) +2 C_0.
$$
Lemma 1 shows that the 
sequence $(a_n)_{n\geq 0}$ is sub-additive: For 
$n,m \geq 0$,
$$
a_{n+m}=d (nx +mx , ny+my)+2C_0\leq d (nx , ny)+d (mx , my)+4 C_0=a_n+a_m.
$$
Thus we have (see lemma 3 below 
for a more general result)
$$
\liminf_{n\to +\infty} {1\over n} a_n =\limsup_{n\to +\infty}
{1\over n} a_n.
$$
Also using $n$ times the inequality from lemma 1 we have
$$
a_n=d (nx, ny)+2C_0\leq n d (x,y)+2(n+1) C_0.
$$
Thus $(a_n/n)_{n\geq 0}$ is a bounded sequence 
and has a limit 
$$
\lim_{n\to +\infty} {1\over n } a_n =\lim_{n\to +\infty} {1\over n }
d (nx, ny) \leq d (x,y)+2C_0.
$$
}

\stit {Remark.}

We need only to use the inequality of lemma 1 
for "large" $x$'s and $y$'s. This will be exploited 
in the proof of theorem 2.

\bigskip

\eno {Lemma 2.}{Let $(E,d)$ be a   
metric vector space over $K$
with the distance $d$ $C_0$-translation 
invariant and $(C_1 ,C_2 , C_3)$-lipschitz
multiplicative.
We define 
$$
d_0 (x,y) =\int_{\dd U} d (u x , u y)\ d \mu (u) .
$$
Then $d_0$ is a distance $(C_1,C_2+C_3)$-lipschitz equivalent 
to $d$, $C_0$-translation 
invariant and $(C_1 ,C_2 , C_3)$-lipschitz
multiplicative. Moreover, $d$ and $d_0$ 
define the same topology on $E$.}

\dem {Proof.}{
Obviously $d_0$ satisfies
the triangle inequality by averaging triangle inequalities.
We  have for $x,y,z \in E$,
$$
\eqalign {
d_0 (x+z, y+z) &\leq \int_{\dd U} d (u (x+z), u (y+z)) \ d\mu (u) \cr
&\leq \int_{\dd U} \left ( d (ux,uy) +C_0\ \right ) d\mu (u) \cr
&\leq d_0 (x,y)+C_0.\cr}
$$
Thus $d_0$ is $C_0$-translation invariant.
Also 
$$
\eqalign {
d_0 (\la x, \la y) &\leq \int_{\dd U} 
d ( u \la x, u \la  y) \ d\mu (u) \cr
&\leq C_1 \ |\la | \ \int_{\dd U} d (ux,uy) \ d\mu (u)
+C_2 |\la |+C_3 \cr
&\leq  C_1 \ |\la |\ d_0 (x,y)+C_2 |\la |+C_3,\cr}
$$
and the reverse inequality follows in the same way. 
Finally the integration over $u\in \dd U$ of
$$
C_1^{-1} d (x,y) -C_2 -C_3 \leq d (ux,uy) \leq 
C_1 d (x,y) +C_2 +C_3,
$$ 
shows  that $d_0$ is $(C_1, C_2+C_3)$-Lipschitz equivalent to $d$.

The distances $d$ and $d_0$ define the same topology. Let 
$(x_n)$ such that $d(x_n , x_0)\to 0$. Then for all $u\in \dd U$
we have $d (ux_n ,ux_0)\to 0$ and the sequence of functions 
$u\mapsto d (u x_n, ux_0)$ are uniformly bounded (the sequence 
$(x_n)$ is $d$-bounded and $\dd U$ is compact). Thus by Lebesgue
dominated convergence we have that $d_0 (x_n , x_0) \to 0$.
Conversely let $(x_n)$ such that $d_0 (x_n , x_0)\to 0$. Then there 
is a sequence $u_n\in \dd U$ such that $d (u_n x_n , u_n x_0)\to 0$.
Since $\dd U$ is compact we can extract a sub-sequence sucht that
$u_n \to u$. Then $d( ux_n , u x_0)\to 0$ so $d (x_n , x_0)\to 0$.
} 

\eno {Definition 7.}{Under the assumptions of theorem 2 
we define 
$$
\d_0 (x,y) =\lim_{n\to +\infty} {1\over n} d_0 (nx, ny)
$$
where $d_0$ is defined in lemma 2,
and 
$$
\d (x,y) = \lim_{n\to +\infty} {1\over n} d (nx, ny).
$$
}

\eno {Lemma 3.}{We have that $\d$ and $\d_0$ are $C_0$-translation 
invariant.}

\dem {Proof.}{It is straightforward from the $C_0$-translation 
invariance of $d$ and $d_0$.}

\eno {Proposition 2.}{We have that $\d$ and $\d_0$ satisfy 
the triangle inequality and are symmetric.
Also if $\lambda \in K$ we have,
$$
\d_0 (\lambda x , \lambda y)=|\lambda| \d_0 (x,y).
$$

}

\eno {Lemma 4.}{Let $\dd Q_+$ be the set of non-negative 
rational numbers. The subset $\dd Q_+ \dd U$ is dense 
in $K$.}

\dem {Proof of Lemma 4.}{Given $\lambda \in K$, $\lambda \not= 0$,
we have $u=\lambda /|\lambda | \in \dd U$. Since the 
restriction of $|.|$ to $\dd Q$ is archimedian, 
we have that $|\dd Q |=|\dd Q_+|$ is dense in $\dd R$,
thus in $|K|$. So there exists a 
sequence of positive rationals $(p_n/q_n)$ such that 
$p_n/q_n  \to |\lambda|$. We conclude that 
$$
\lim_{n\to +\infty} {p_n \over q_n} u =\lambda.
$$
}

\dem {Proof of proposition 2.}{The triangle inequality and the 
symmetry is immediate from the definition.

Given an integer $p\geq 0$ we have
$$
\d_0 (px,py) =\lim_{n\to +\infty} {1\over n} d_0 (npx,npy)
=p \lim_{n\to +\infty} {1\over np} d_0 (npx,npy)=p \d_0 (x,y).
$$

Now if $p/q\in \dd Q$, $q\geq 1$, $p\geq 1$ we have
$$
q \d_0 (p/q \  x, p/q \ y)=p q \d_0 (1/q\ x , 1/q \ y)
=p \ \d_0 (x,y).
$$
so for any rational number $r\in \dd Q_+$,
$$
\d_0 (rx , ry) =|r|\ \d_0 (x, y).
$$
}

\eno {Proposition 3.}{We have 
$$
\d_0 (x,y) =\int_{\dd U } \ \d (ux, uy) \ d\mu (u).
$$
}

\dem {Proof of proposition 3.}{From Lemma 1 and 
compactness of $\dd U$ we have 
that 
$$
{1\over n } d (nux, nuy) \leq d (ux,uy) +2 C_0 \leq M (x,y),
$$
where $ M(x,y)$ is a bound uniform on $n$. Thus the functions
$u\mapsto {1\over n} d (nux, nuy)$ are uniformly bounded. 
From Lebesgue's dominated convergence
theorem we have
$$
\eqalign {
\d_0 (x,y) &=\lim_{n\to +\infty} {1\over n} d_0 (nx,ny) \cr
&=\lim_{n\to +\infty} \int_{\dd U} d (unx, uny) \ d\mu (u) \cr
&=\lim_{n\to +\infty} \int_{\dd U} d (nux, nuy) \ d\mu (u) \cr
&=\int_{\dd U} \lim_{n\to +\infty} d (nux, nuy) \ d\mu (u) \cr
&=\int_{\dd U} \d (ux , uy) \ d\mu (u)\cr
}
$$
q.e.d.}

\eno {Lemma 5.}{Any closed (resp. open) set for $\d_0$
 is closed (resp. open) set for $d.$ So the topology generated by $d$
is richer than the topology generated by $\d_0.$}

\dem {Proof of Lemma 5.}{Let $x_n$ be a sequence 
of points in $E$ such that 
$d(x_n,x)\to 0$.

We have 
$$ 0 \leq \d_0(x_n,x) \leq C_1 d(x_n,x) +C_2.$$
Therefore
$$ \limsup_{n\to +\infty} \d_0(x_n,x) \leq C_2.$$

By proposition 2, for all $ \lambda \in K$,
$$\d_0(x_n,x)={1 \over |\lambda|}\d_0(|\lambda|x_n,|\lambda|x)
\leq {1 \over |\lambda|}(C_1 \d_0(|\lambda|x_n,|\lambda|x)+C_2).$$
Thus
 $$0 \leq \limsup_{n\to +\infty} \d_0(x_n,x) \leq
{ C_2 \over |\lambda|}$$
for all $\lambda \in K$. Taking limit when $|\lambda|$ tends to $+\infty$
we obtain 
$$0 \leq \limsup_{n\to +\infty} \d_0(x_n,x) \leq 0.$$
Thus $\lim_{n\to +\infty} \d_0(x_n,x)=0$ and the lemma follows.}

\eno {Proposition 4.}{Let
 
\item {(i)} $E_0 =\{ x\in E ; {\hbox {\rm for all }} u\in \dd U, 
\d (ux,0) =0 \}$,

\item {(ii)} $E_1 =\{ x \in E ; \d_0 (x,0) =0 \}$,

\item {(iii)} $E_2$ maximal subspace where $d_0$ is bounded,

\item {(iv)} $E_3$ maximal subspace where $d$ is bounded,

Then $E_0=E_1=E_2=E_3$.}

\dem {Proof of proposition 4.}{We have 
$$
\d_0 (x,0) =\int_{\dd U} \d (ux,0) \ d\mu (u),
$$
thus  $E_0 = E_1$.
Observe that $E_1$ is a subspace of $E$, for $x,y \in E_1$,
using the translation invariance,
$$
\d_0 (x+y, 0) \leq \d_0 (x,0) +\d_0 (y,0) =0+0=0.
$$
Also $\d_0 (\lambda x, 0)=\lambda \d_0 (x,0) =0$.
Moreover for $x,y\in E_1$ we have
$$
d_0  (x,y) \leq C_1 \left ( \d_0 (x,y) +C_2  \right )
\leq  C_1 \left ( \d_0 (x,0) +\d_0 (y,0) +C_2 \right )\leq C_1 C_2 .
$$
Thus $d_0$ is bounded in $E_1$ and $E_1\subset E_2$. 
From the definition of $\d_0$ it follows that $E_2 \subset E_1$,
thus $E_1=E_2$. 

Finally $E_2=E_3$ because $d_0$ and $d$ are lipschitz equivalent.}

\eno {Proposition 5.}{The subspace $E_0$ is a closed 
subspace of $E$.}

\dem {Proof of proposition 5.}{let $(x_n)$ be a converging sequence of points
in $E_0$, $x_n\to x$. If $x\notin E_0$ then $\d_0 (x,0)\not=0$. Let 
$C=C_1 C_2 +C_1 C_3 +C_2$. Consider 
$$
y_n ={C +1 \over \d_0 (x,0)} x_n.
$$
We have $y_n \in E_0$ (since $E_0$ is a subspace), $y_n\to y$ with 
$$
\d_0 (y ,0) ={C +1 \over \d_0 (x,0)} \d_0 (x,0)> C.
$$
Then
$$
\eqalign {
C < \d_0 (y,0)&\leq \d_0 (y-y_n , 0) +\d_0 (y_n ,0)=\d_0 (y-y_n ,0) \cr
&\leq C_1 d_0 (y, y_n )+C_2 \leq C_1^2 \  d (y,y_n) +C_1 (C_2 +C_3) +C_2=
C_1^2 \  d (y,y_n) +C. \cr}
$$
Passing to the limit $n\to +\infty$, we get $C<C$. Contradiction.}  

\eno {Proposition 6.}{We define for a class $\bar x \in E/E_0$,
$$
||\bar x ||=\d_0 (x,0).
$$
The definition is independent of the representant $x$ of the 
class $\bar x =x+E_0$ and $||.||:E/E_0 \to \dd R_+$ is a norm.}

\dem {Proof of proposition 6.}{If $y\in \bar x$ then $x-y\in E_0$ thus
$\d_0 (x-y, 0)=0$ so by translation invariance $\d_0 (x,y)=0$
and 
$$
\d_0 (x,0)\leq \d_0 (x,y) +\d_0 (y,0)=\d_0 (y,0).
$$
In the same way $\d_0 (y,0)\leq \d_0 (x,0)$ and finally
$\d_0 (x,0)=\d_0 (y,0)$.
Also if $||\bar x ||=0$ then $\d_0 (x,0)=0$ and $x\in E_0$, i.e.
$\bar x =\bar 0$. The other properties of a norm follow 
from the properties of $\d_0$.}

\eno {Definition 8.}{We denote by $D : E/E_0 \to \dd R_+$, resp.
$D_0 : E/E_0 \to \dd R_+$, the Hausdorff distances for $d$, resp. $d_0$, 
between classes modulo $E_0$,
$$
D (x_0+E_0 , y_0+E_0) =\max \left ( \inf_{x\in \bar x_0} 
\sup_{y\in \bar y_0} d (x,y) , \inf_{y\in \bar y_0} 
\sup_{x\in \bar x_0} d (x,y) \right ).
$$ 
}

Hausdorff distances over non-compact sets do not need to be 
finite. Also when the distance is zero the sets do not need
in general to coincide. In our situation Hausdorff
distances define a proper distance on classes  modulo $E_0$.

\eno {Lemma 6.}{The spaces $(E/E_0 , D)$ and $(E/E_0, D_0)$ are
metric spaces, and $D$ and $D_0$ define the quotient
topology on $E/E_0$.}
\dem {Proof of lemma 6.}{We first prove that $D$ and $D_0$ do define distances.
We carry the proof for $D$. The same proof applies 
to $D_0$. If $x\in \bar x_0$ and $y\in \bar y_0$ we have
$$
d (x,y)\leq d(x,x_0) +d(x_0 ,y_0) +d(y_0,y)
\leq d (x_0 , y_0) +d (x-x_0, 0) +d (y_0-y,0)+ 2C_0,
$$
and the last two terms are uniformly bounded since $x-x_0 \in E_0$
and $y_0-y \in E_0$. This shows that $D (\bar x_0 , \bar y_0) <+\infty$.

Assume that $D(\bar x_0 , \bar y_0)=0$. Then there is a sequence
$(e_n)$ with $e_n \in E_0$ such that $y_0+e_n \to x_0$. Therefore
$e_n \to x_0-y_0$. But we have proved that $E_0$ is closed,
thus $x_0-y_0 \in E_0$ and $\bar x_0 =\bar y_0$.

We denote $\pi : E \to E/E_0$ the quotient map.
We prove that each open set $U$ for the quotient topology of
$E/E_0$ is open for $D$. Let $\bar x_0 \in U$ and $D (\bar x_n ,
\bar x_0 )\to 0$. We have to prove that there exists $N$ such 
that for $n\geq N$, $\bar x_n \in U$. We have that $x_n +E_0 
\to x_0 +E_0$ in Hausdorff metric for $d$. Thus there exists
a sequence $x_n'\in x_n +E_0$ such that $d(x_n', x_0)\to 0$.
Since $d$ defines the topology of $E$ and $\pi^{-1} (U)$ is open
there exists $N$ such that for $n\geq N$ we have $x_n' \in 
\pi^{-1} (U)$. Then $\bar x_n =\bar x_n' =\pi (x_n')\in U$.q.e.d.
}

\eno {Proposition 7.}{The norm $||.||$ is $(C_1 ,C_2+C_3)$-lipschitz
equivalent to $D_0$ and $(C_1^2 , C_2')$-lipschitz equivalent 
to $D$.}

\dem {Proof of proposition 7.}{We have for $x\in \bar x_0$, $y \in \bar y_0$,
$$
C_1^{-1} d_0 (x,y) -(C_2 +C_3) \leq \d_0 (x,y)=\d_0 (x_0, y_0) =
||\bar x_0 -\bar y_0 ||
\leq C_1 d_0 (x,y)
+(C_2 +C_3 ).
$$
Now letting $x$ and $y$ run over $\bar x_0$ and $\bar y_0$ 
respectively, and using the definition of $D$ we have the 
result. Same proof for $D_0$.} 
 
\null
\vfill
\eject

\tit {2) Proof of theorem 3.}

\stit {The conditions are necessary.}

We assume that $d$ is asymptotically isometric to a 
norm $||.||$. Let $C_1 >1$. Then $C_1^{1/2} >1$ and 
there exists $C_0 \geq 0$
such that for $x,y\in E$,
$$
C_1^{-1/2} ||x-y||-C_0 \leq d (x,y) \leq C_1^{1/2} ||x-y|| +C_0.
$$
Then for any $\la \in K$ we have,
$$
\eqalign {
d (\la x , \la y ) \leq C_1^{1/2} ||\la x -\la y||+C_0
&= C_1^{1/2} |\la |\ || x - y||+C_0 \cr
&\leq C_1^{1/2} |\la | \left ( C_1^{1/2} d (x,y) 
+C_1^{1/2} C_0 \right ) +C_0 \cr
&\leq C_1 |\la |\ d (x,y) +C_1 C_0 |\la | +C_0 \cr
}
$$
The reverse inequality is proved in the same way 
and $d$ is asymptotically multiplicative.
Also for $x_1, \ldots x_n , y_1, \ldots , y_n \in E$
we have 
$$
\eqalign {
d \left ( \sum_{i=1}^n x_i , \sum_{i=1}^n y_i \right )
&\leq C_1^{1/2} || \sum_{i=1}^n x_i -\sum_{i=1}^n y_i ||+C_0 \cr
&\leq C_1^{1/2} \sum_{i=1}^n ||x_i-y_i|| +C_0 \cr
&\leq C_1 \sum_{i=1}^n d (x_i , y_i) + nC_1 C_0 +C_0\cr
&\leq C_1 \sum_{i=1}^n d (x_i , y_i) + n (2 C_1 C_0)  \cr
}
$$
thus the condition in theorem 2 is necessary.

\stit {The conditions are sufficient.}

We construct the norm by the same strategy as in theorem 1.
We need a refinement on the lemma on sub-additive
sequences.

\eno {Lemma 7.} {Let $(a_n)_{n\geq 0}$ be a sequence
of real numbers satisfying 
the following weak sub-additive property: For any $C_1>1$ 
there exists $C_0\geq 0$ such that for any $q\geq 2$ and any 
$m_1, \ldots , m_q \geq 0$,
$$
a_{m_1 +\ldots +m_q} \leq C_1 \sum_{i=1}^q a_{m_i} +q C_0.
$$
Then 
$$
\liminf_{n\to +\infty} {1\over n} a_n =
\limsup_{n\to +\infty} {1\over n} a_n .
$$
}
\dem {Proof.}{Fix for the moment $C_1 >1$. 
Choose  $n\geq 1$. For any $m\geq 0$ we can consider
the euclidian division $m=nq+r$, with $0\leq r <n$.
We have  
$$
a_m=a_{nq+r} \leq C_1 (qa_n +a_r )+ (q+1) C_0.
$$
Dividing by $m$ and taking the least upper bound for 
$m\to + \infty$ ($q\to + \infty$) we have
$$
\limsup_{m\to +\infty } {1\over m} a_m 
\leq C_1 {1\over n } a_n  + {C_0 \over n} .
$$
Now taking the greater lower bound for $n\to +\infty$
we get
$$
\limsup_{n\to +\infty } {1\over n} a_n
\leq C_1 \liminf_ {n\to +\infty } {1\over n} a_n.
$$
Since this holds for any $C_1>1$ the lemma follows.}

It is simple to check that we have the same lemma as 
lemma 2 for theorem 1:

\eno {Lemma 8.}{We assume the hypothesis of theorem 2.
We define 
$$
d_0 (x,y) =\int_{\dd U} d (ux,uy) \ d\mu (u) \ .
$$
Then $d_0$ is a distance assymptotically equivalent to $d$ 
and satisfying the hypothesis of theorem 3.}

Now we have:
 
\eno {Proposition 8.}{For any $x,y\in E$, the limit 
$$
\d (x,y) =\lim_{n\to +\infty} {1\over n} d_0 (nx, ny)
$$
exists.}

\dem {Proof.}{The sequence 
$$
a_n =d_0 (nx, ny)
$$
is weakly sub-additive:

$$ \eqalign {
a_{m_1+ \cdots +m_q} =
& d_0 \left ( (m_1+ \cdots +m_q)x, \ (m_1+ \cdots +m_q) y \right ) \cr
&= d_0 (m_1 x + \cdots +m_q x , \ m_1 y + \cdots +m_q y) \cr
& \leq C_1 \sum_{i=1} ^q a_{m_i} +qC_0. \cr
}$$

 Moreover $a_n/n$ is bounded:
$$ \eqalign {
a_n = d (nx, \ ny) 
&= d (x+ \cdots +x, \ y+ \cdots +y) \cr
&\leq C_1 \sum_{i=1} ^n d (x, \ y) + nC_0 \cr
& =  C_1 n d (x, \ y)+nC_0. \cr
}$$
 Thus 
$${a_n \over n} \leq C_1 d (x, \ y)+C_0.$$

 The result 
follows from lemma 7.}

\eno {Proposition 9.}{For $x,y\in E$ we define 
$$
\d_0 (x,y)=\lim_{n\to +\infty} {1\over n} d_0 (nx,ny).
$$
Then $\d_0$ is translation invariant. If we define
$$
||x-y||=\d_0 (x,y)
$$
then $||.||$ is a norm that is asymptotically isometric 
to $d_0$, so also to $d$.}

\dem {Proof.}{We prove the translation invariance. The rest
follows the same lines as the proof of proposition 2.
Let $x,y,z \in E$. For any $C_1 >1$, there exists $C_0 \geq 0$
such that 
$$  \eqalign {
d_0(n(x+z), n(y+z))
& =d_0 (nx+nz, ny+nz) \cr
& \leq C_1 d_0(nx,ny)+C_1 d_0 (nz,nz) +2C_0 \cr
& =C_1 d(nx,ny)+2C_0 \ . \cr
}
$$
Now dividing by $n$ and passing to the limit $n\to + \infty$, 
we get, for any $C_1>1$,
$$
\d_0 (x+z, y+z) \leq C_1 \d_0 (x,y) \ .
$$
Therefore making $C_1 \to 1$, for all $x,y,z\in E$,
$$
\d_0 (x+z, y+z)\leq \d_0 (x,y) \ .
$$
Replacing $x$ by $x+z$, $y$ by $y+z$ and $z$ by $-z$, we get
the opposite inequality, and the translation invariance.}

This finishes the proof of theorem 3.
\bigskip

\tit {Theorem 3 implies Theorem 1.}

Let $(E,d)$
be asymptotically isometric to an AMN vector space. Therefore
$d$ is unbounded in non-trivial subspaces and we have seen
that it satisfies the other
hypothesis of theorem 3. Thus $d$ is 
asymptotically isometric to a norm of $E$.

\null
\vfill
\eject

\tit {\centerline {\bf REFERENCES}} 

\bigskip

\ref{[GH]}{E. GHYS, P. DE LA HARPE}{Sur les groupes hyperboliques
d'apr\`es M. Gromov}{Progress in Mathematics, {\bf 83}, 
Birkh\"auser, Boston, 1990.}

\ref{[Gr]}{M. GROMOV}{Hyperbolic groups}{"Essays in group theory",
G. Gernsten editor, Mth. Sci. Res. Inst. Pub., Springer, 1987, p. 75-263.}

\ref{[Mu]}{E. MU\~NOZ GARCIA}{H\"older absolute values are 
equivalent to classical ones}
{Proc. Amer. Math. Soc.  {\bf 127}, 1999, p. 1967-1971.}

\ref{[Sch]}{H.H. SCHAEFER}{Topological Vector Spaces}
{Graduate texts in Mathematics, {\bf 3}, Springer-Verlag, 1970.}

\ref {[We]}{A. WEIL}{Introduction to number theory}
{Springer-Verlag.}

\end